\documentclass[10pt]{amsart}
\usepackage[leqno]{amsmath}
\usepackage{hyperref}
\usepackage{appendix}
\usepackage{graphicx}
\usepackage{mathrsfs}
\DeclareGraphicsExtensions{.pdf,.png,.jpg}

\newcommand{\tuple}[1]{\ensuremath{\left \langle #1 \right \rangle }}

\theoremstyle{plain}
  \newtheorem{theorem}{Theorem}

  \newtheorem*{churchthesis}{Church's Thesis}
  \newtheorem*{turingthesis}{Turing's Thesis}
  
  \newtheorem*{argInf}{Argument by Infinity}
  \newtheorem*{argVague}{Argument by Physical Vagueness}
  \newtheorem*{thesisP}{Thesis P}
  \newtheorem*{thesisM}{Thesis M}
  \newtheorem*{use-constr}{Usability Constraint}
  \newtheorem*{superT}{Super-Turing Computation}
  \newtheorem*{pt}{Post-Turing Thesis}

\theoremstyle{definition}

\theoremstyle{remark}

\numberwithin{equation}{subsection}

\title{Practical intractability: a critique of the hypercomputation movement}
\author{Aran Nayebi}
\email{anayebi@stanford.edu}
\urladdr{http://www.stanford.edu/~anayebi}
\begin{document}
\begin{abstract}
For over a decade, the hypercomputation movement has produced computational models that in theory solve the algorithmically unsolvable, but they are not physically realizable according to currently accepted physical theories. While opponents to the hypercomputation movement provide arguments against the physical realizability of specific models in order to demonstrate this, these arguments lack the generality to be a satisfactory justification against the construction of \emph{any} information-processing machine that computes beyond the universal Turing machine. To this end, I present a more mathematically concrete challenge to hypercomputability, and will show that one is immediately led into physical impossibilities, thereby demonstrating the infeasibility of hypercomputers more generally. This gives impetus to propose and justify a more plausible starting point for an extension to the classical paradigm that is physically possible, at least in principle. Instead of attempting to rely on infinities such as idealized limits of infinite time or numerical precision, or some other physically unattainable source, one should focus on extending the classical paradigm to better encapsulate modern computational problems that are not well-expressed/modeled by the closed-system paradigm of the Turing machine. I present the first steps toward this goal by considering contemporary computational problems dealing with intractability and issues surrounding cyber-physical systems, and argue that a reasonable extension to the classical paradigm should focus on these issues in order to be practically viable.
\end{abstract}

\maketitle
\tableofcontents
\newpage
\section{Introduction}
\indent The seminal work of G\"{o}del, Church, Turing, and Post in the 1930s regarding the analysis of algorithmic process provided a framework for future theoretical and practical investigations in computer science. Letting ``effective calculability'' refer to an informal notion, and ``computable'' to mean ``computable by a Turing machine''\footnote{In places where clarification may be necessary, I will use the equivalent term ``Turing-computable''. When referring to the non-Turing computable, I will simply use the standard term ``non-computable''.}, Church's \cite{church} and Turing's \cite{turing} analyses of effective calculability ultimately led to
\begin{churchthesis}
A function of positive integers is effectively calculable iff it is recursive.
\end{churchthesis}
\begin{turingthesis}
What is effectively calculable is computable.
\end{turingthesis}
When restricted to functions of positive integers, both Church's thesis and Turing's thesis are equivalent notions, hence the term ``Church-Turing thesis'' which was introduced by Kleene \cite{kleene1967} in 1967. Kleene \cite{kleene1952} provides extensive evidence in support of the Church-Turing thesis -- one important justification being that every known effectively calculable function has been shown to be computable. Turing had in mind calculation by an idealized human being who proceeds mechanically with pencil and paper. A machine in the sense of Turing is essentially a finite state device operating on a one-way finite (but potentially unbounded) linear tape divided into squares, meant to be deterministic and sequential, comparable to ``a man in the process of computing a real number [in that it is] only capable of a finite number of conditions'' \cite[pg. 231]{turing}. Although the Turing machine was meant as a \emph{gedankenexperiment} for proving the algorithmic unsolvability of the \emph{Entscheidungsproblem}, the early computers of the 1950s were, in some sense, physical realizations of this idealization. Thus, the algorithmically unsolvable problems of the theoretical realm were generally believed to remain unsolvable in the practical realm as well, namely, one could not build a physical device that could solve the algorithmically unsolvable. Though theoretical devices that can compute more than the universal Turing machine have been around since 1939 (starting with Turing's oracle machine \cite{turing1939}), the ``hypercomputation'' movement has been relatively recent. \newline
\begin{quote}
A hypercomputer is any information-processing machine, notional or real, that is able to achieve more than the traditional human clerk working by rote. Hypercomputers compute functions or numbers, or more generally solve problems or carry out tasks, that lie beyond the reach of the universal Turing machine of 1936. Copeland \cite[pg. 462]{copeland2002} \newline
\end{quote}
\indent One aim of hypercomputation has been to argue for the existence of \emph{physical} devices that can compute non-computable functions. The devices presented so far do not work in practice because these hypercomputational models rely on embedded infinities to perform computations, or on some black box oracle that computes a non-recursive function, without much (or any) description of how this oracle explicitly goes about its computations. It is common for opponents of the hypercomputation movement to point to embedded infinities in order to argue against the physical realizability of \emph{specific} models, but the issue with hypercomputation is far deeper and more general than that. \newline
\indent I will argue that the aims of the hypercomputation movement are not well-defined. Thus, what is needed is a deeper analysis of Turing, and I will argue that any sort of ``plausible hypercomputability'' would have to explain what it is that it expects to compute that is outside the scope of Turing \cite{turing} in at least as mathematically precise a manner as Turing did. I will demonstrate that attempts to make the current aims of the hypercomputation movement mathematically precise immediately lead to physical impossibilities given current physical theories, and this gives us impetus to shift the current focus of hypercomputation. Instead of asserting that non-computable functions, tasks, etc. are physically computable, a more reasonable starting point for an extension to the classical paradigm should be focused on modeling processes that are not adequately encapsulated by the Turing machine as a closed system. As guiding examples, I consider contemporary computational problems dealing with issues surrounding cyber-physical systems as well as intractability. The question of whether or not these processes allow one to physically compute the non-computable is not considered.

\section{Models of hypercomputation}
\indent Can one, given current physical theories, build a machine that solves, say, the halting problem or some other function not computable by the universal Turing machine? Since there have been many such proposed models (and at least speculation of physical  processes that cannot be simulated by Turing machine has spanned over four decades\footnote{Copeland \cite[pp. 130-1]{cope98} provides an excellent list of papers concerned that present these speculations prior to the coining of the term ``hypercomputation'' by Copeland and Proudfoot \cite{copeland-proudfoot} in 1999.}), the latter question can be reformulated as an inquiry into the physical realizability of these models. The most well-known critique of hypercomputation and negative answer to the above question is provided by Martin Davis \cite{davis2004}\cite{davis-sr}\cite{davis2006}. His arguments (and those of other authors) against these models primarily fall under what I will call the ``Argument by Infinity'' and the ``Argument by Physical Vagueness''.
\begin{argInf}
A model that involves embedded infinities is not physically realizable.
\end{argInf}
\begin{argVague}
It is an open question/unclear how to take advantage of certain physical phenomena required by a particular model in order to compute the non-computable $($not to mention that there may be questions about the existence of the physical phenomena$)$, or how to construct the particular hypercomputer based on the physical theory $($or theories$)$ that it is claimed to be consistent with.
\end{argVague}
I refrain from surveying various models of hypercomputation, as many surveys on the subject can be found, such as by Copeland \cite{copeland2002} or by Ord \cite{ord}. Note that it is impractical to categorize all arguments against the existence of hypercomputers due to the rather large literature on the subject, but these two cover most (if not all) of them. Of course, the Arguments by Infinity and Physical Vagueness do not just apply to models of hypercomputation but to \emph{any} computational model that involves infinities or is physically vague. The Argument by Physical Vagueness is quite often the consequence of a model's attempt to utilize an infinite process to compute the non-computable; thus, it is not unusual that the two arguments go together. Furthermore, it should be noted that the Argument by Physical Vagueness \emph{alone} is weaker than the Argument by Infinity, because unlike the latter, the former does not rule out a model (since we are operating under the plain assumption that the physically constructible is finite\footnote{Though one can never exclude the possibility that a groundbreaking new physical theory could provide for the physical implementation of the infinite; however, so far no such theory is even closely there.}), but just makes its existence in the near future highly unlikely. \newline
\indent As one can see, the Argument by Infinity and the Argument by Physical Vagueness are \emph{model-specific}, namely, one applies it to a particular model to rule out its physical plausibility. They by no means provide a satisfactory argument against the physical construction of \emph{any} information processing machine that computes the non-computable given current physical theories. However, the fact that hypercomputational models presented so far have not been shown to be physically feasible due to the Argument by Infinity and the Argument by Physical Vagueness gives one motivation to pursue more general arguments against hypercomputation, that I will do in \S 3. \newline

\section{A general negative answer to physical hypercomputation}
\indent One cannot point to a clear answer to the question, ``What is hypercomputation?''. The quote from Copeland presented in \S 1 is exceptionally vague because when one attempts to clarify what is meant by an ``information processing machine'' that computes functions or solves problems that ``lie beyond the reach of the universal Turing machine'', one is only faced with a plethora of models that supposedly proceed to do this in a variety of different manners, as evidenced in the surveys of \cite{copeland2002} or \cite{ord}, and even the models explored in \S 3.2. Moreover, the common theme between these models is that they are often based on the trivial insight that if a physical process exists that generates a non-computable function (or say, an infinite precision real number), then it could be exploited (somehow) to compute the non-computable. And despite this triviality, there is the further catch that in exploiting such physical processes one falls victim primarily to the Argument by Infinity and/or the Argument by Physical Vagueness.\newline
\indent However, as Stannett comments, ``those in favour of hypercomputation can sometimes be too quick to assume that their models do not include hidden infinities, while those against hypercomputation are sometimes too quick to assume that what they consider to be the ingredients of computation are in some sense `necessary' ''\cite[pg. 9]{stannett}. Thus, what is developed here is a more mathematically concrete challenge to hypercomputability that is not model-specific, and I will show that one is immediately led into physical impossibilities, therefore providing an effective starting point toward demonstrating the physical impracticality of hypercomputers more generally. \newline
\indent The best place to start is with a deeper analysis of Turing \cite{turing}. Since hypercomputation aims to compute that which lies outside the scope of the universal Turing machine, then what is needed is a mathematically clear conception of that scope. \newline
\subsection{A deeper analysis of Turing's and Gandy's arguments}
\indent Turing's analysis in \cite{turing} is essentially a theory of \emph{human} computation. In \S 1 of his paper, Turing briefly describes what he means by ``computable real number'' as one whose decimal expression is ``calculable by finite means'' \cite[pg. 230]{turing}. This definition as well as the comparison in \S 1 to a man in the process of computing a real number to the computation performed by an $a$-machine (which is what has become known as the ``Turing machine''\footnote{It was first called this by Church in his 1937 review \cite{church1937} of Turing's paper.}) is not justified until \S 9 of his paper. \newline
\indent In this justification, one finds that he provides three kinds of arguments: ``(a) A direct appeal to intuition. (b) A proof of the equivalence of two definitions (in case the new definition has a greater intuitive appeal). (c) Giving examples of large classes of numbers which are computable'' \cite[pg. 249]{turing}. It is evident then that Turing's emphasis when describing the formal mathematical notion of computability is on how much it ``agrees'' with one's intuition of effective calculability, where ``effective'' is taken to mean that the calculation is deterministic and terminates after a finite amount of time. In fact, this direct appeal to intuition (the first argument) is what is most remembered, in spite of the fact that Turing gave two other arguments in support of his justification, as mentioned above. This intuitive definition is widely accepted, one important reason being its convincing psychological appeal. For instance, Church \cite{church1937} found Turing's method regarding the unsolvability of the \emph{Entscheidungsproblem} to be more favorable than his own. Moreover, G\"{o}del praised Turing's explication stating that it is ``absolutely impossible that anybody who understands the question and knows Turing's definition should decide for a different concept'' \cite[pg. 84]{wang1974}. \newline
\indent Despite this praise, an appeal to intuition is not as mathematically rigorous as possible. Instead, what is needed is an axiomatization of Turing's explication to serve as guiding principles for the inquiry of this section. Of course, this axiomatization would be implicit in Turing's original argument, even though he did not explicitly state it. Certainly Turing's focus on \emph{human} mechanical calculations was appropriate given the historical context of his work, for any solution of the \emph{Entscheidungsproblem} was to be carried out by us! The question that immediately arises, however, is that even though Turing started by analyzing human computation, why should one follow him in that? The reason is that Turing's analysis offers a detailed and precise investigation that connected an informally understood concept (that of ``effective calculability'') to a mathematically precise notion (that of ``computability''). And this analysis, through its use of restrictions, allows us to readily segue into the notion of \emph{machine computation}, which is of prime relevance to the investigative aims of this section regarding a more mathematically precise understanding of hypercomputability. \newline
\indent One primary reason why Turing's analysis of effective calculability made the previously informally understood concept mathematically rigorous was that he placed constraints on the steps permitted in calculations -- these were entirely based on the limitations of the human computing agent. Sieg \cite{sieg2008} formalizes these limitations as such: 
\begin{enumerate}
\item \textbf{Boundedness:} A computor can immediately recognize only a bounded number of configurations.
\item \textbf{Locality:} A computor can change only immediately recognizable configurations,
\end{enumerate}
where, under Sieg's terminology, a ``computor'' refers to a human computing agent who proceeds mechanically. It is easy to see that processes satisfying the above conditions are computable by Turing machines since the first conditions require that the configurations, on which the operations are performed, have a fixed bounded size; and the second conditions require that the operations take place locally, within a fixed finite neighborhood. The axioms that Sieg \cite{sieg2008} derives for Turing computors are motivated by these boundedness and locality conditions; I will not state them here.
\newline
\indent It is Gandy \cite{gandy} who makes explicit the connection between the abstract human being working in a rote manner to the physical construction of a machine. While Turing viewed the aforementioned restrictive provisos due to the memory limitations of computors (stating that the ``human memory is necessarily limited'' \cite[pg. 231]{turing}), Gandy, on the other hand, appealed to evident physical limitations of machines as impetus to derive restrictive principles for discrete deterministic mechanical devices. By considering and analyzing cellular automata and parallel computation, he produces four principles that any (non-analog) machine must satisfy and proves that devices satisfying them are Turing-equivalent.\newline
\indent The restrictive conditions that Gandy presents are motivated by the physical considerations of quantum mechanics and special relativity, namely, ``that there is a lower bound on the linear dimensions of every atomic part of the device and that there is an upper bound (the velocity of light) on the speed of propagation of changes'' \cite[pg. 126]{gandy}. Gandy aims to provide an argument for the following thesis:
\begin{thesisM}
What can be calculated by a machine is computable.
\end{thesisM}
He clarifies that he will consider only discrete mechanical devices, and proceeds to formulate the notion of a discrete mechanical device in terms of four precise axiomatic principles. Devices (synonymous with ``machine'') satisfying Gandy's principles (termed ``Gandy machines'') are based on a set-theoretic framework for specification, consisting of classes of ``states'' with a ``transition operation'' allowing for a step-wise transition from one state to the next. States are represented by subclasses of hereditarily finite sets ($HF$) built up from a potentially infinite set of atoms, conceived as the basic components of the machine, closed under isomorphic structures (such subclasses are called ``structural classes''). Transitions are given by restricted operations from states to states, termed ``structural operations''. Ignoring the technical details of Gandy's presentation (which are clarified in Sieg and Byrnes \cite{sb1999b} and Sieg \cite{sieg2002}), one can manage to approximate the four principles as follows: 
\begin{enumerate}
\item \textbf{Principle I (The Form of Description):} A Gandy machine $M$ can be described by the pair $\tuple{S, F}$, where $S$ is a structural class, and $F$ is a structural operation from $S$ to $S$. If $x_0 \in S$ is an initial state of $M$, then $F(x_0)$, $F(F(x_0))$, $\ldots$ are its subsequent states. \newline \newline
The remaining three principles place restrictions on $S$ and $F$ involving certain finiteness or boundedness conditions, and Gandy \cite{gandy} shows that if any of these principles are significantly weakened in almost any way, then every number-theoretic function becomes calculable.
\item \textbf{Principle II (The Principle of Limitation of Hierarchy):} The set-theoretic rank of the states is bounded. Roughly speaking, each state of $S$ can be assembled from parts, and those parts can be assembled from other parts, etc., but the complexity of this structure is bounded.
\item \textbf{Principle III (The Principle of Unique Reassembly):} Each state of $S$ can be assembled from basic parts of bounded size, and there is a unique way of putting these parts together (while model construction-kits aim to satisfy this principle, this principle also takes into consideration parts which do overlap).
\item \textbf{Principle IV (The Principle of Local Causality):} For any state $x \in S$, $F(x)$ is assembled from bounded (and possibly overlapping) parts of $x$. This is an abstraction of the aforementioned physical considerations of quantum mechanics and special relativity that Gandy takes into account.
\end{enumerate}
From these four principles, Gandy \cite{gandy} proves the following
\begin{theorem}
What can be calculated by a device satisfying Principles I-IV is computable.
\end{theorem}
As a result, he clarifies his arguments for a more definite version of ``Thesis M'' in the following way:
\begin{thesisP}
A discrete deterministic mechanical device satisfies Principles I-IV.
\end{thesisP}
\indent In sum, \emph{Thesis P} (or \emph{Gandy's thesis}) is the assertion that any discrete deterministic mechnical device must satisfy the restrictive principles of local causation and unique assembly, and as Sieg and Byrnes \cite{sb1999b} point out, this is analogous to Turing's Thesis which claims that any human computor, which can be represented by a discrete dynamical system and his operation, must obey appropriate locality conditions. Of course, the important difference between the Gandy machine and the Turing machine is that latter can modify only a bounded part of a state, whereas the former allows state-transitions that result in arbitrarily many bounded parts, since Gandy machines take parallel working into account. \newline
\indent Despite the shortcoming that Principle IV does not apply to machines obeying Newtonian mechanics, since ``in these there may be rigid rods of arbitrary lengths and messengers travelling with arbitrary large velocities, so that the distance they can travel in a single step is unbounded'' \cite[pg. 145]{gandy}, Gandy's principles are a reasonable (in the sense that they obey the modern physical theories of quantum mechanics and special relativity) and mathematically concrete formulation of physical principles that a real machine\footnote{I am only interested in the consideration of real machines here, though Gandy's Thesis leaves open to interpretation the possible consideration of both notional and real machines (cf. Shagrir \cite[\S 4]{humach}).} must satisfy. Not to mention that since the early twentieth century, it has been known that at extremes of velocity (for instance, near the speed of light) and on incredibly small scales, Newtonian mechanics is not a successful physical theory. Moreover, it is worth mentioning that Sieg \cite{sieg2008} obtains an axiomatization of the class of Gandy machines based on Gandy's principles of local causation and unique assembly, and Sieg's axioms for discrete dynamical systems (for both human and machine computations) result in the reduction of models of these axioms to Turing machines, and cellular automata and many artificial neural nets, for example, can be shown to satisfy these axioms \cite{sieg2008}. This axiomatization thus provides good evidence for Gandy's principles for mechanisms as sufficiently general, and in fact, in \cite{sb1999b} and \cite{sieg2002}, two formulations of abstract models for parallel computations are given based on Gandy's \cite{gandy} precise mathematical framework (I will mention their relation to other paradigms of computing in \S 4.2.1). 
\subsection{Addressing arguments against Gandy's thesis}
Copeland and Shagrir \cite{copeland-shagrir}, however, argue that Gandy's principles for mechanisms are not general enough, and instead they present ``(ideal) physical machines that lie outside the class of `Gandy machines' and that compute (in a broad sense of  `compute') non-recursive functions....[including] discrete, deterministic hypercomputers from both Newtonian and relativistic frameworks, as well as non-discrete and non-deterministic hypercomputers'' \cite[pg. 220]{copeland-shagrir}. The focus of this paper is, of course, on the physical possibility of hypercomputation, not the notional one. Although the question of whether discrete deterministic mechanical devices contravening Gandy's principles are allowed by the physics of the real world is a question that is empirical not logical, which the authors are apt to point out, it is precisely the logical framework of Gandy's that makes the empirical affirmation of hypercomputation unlikely, as I will show below. \newline
\indent The machines that Copeland and Shagrir present that lie outside the class of Gandy machines comprise none other than a partial survey of hypercomputational models in the literature. Before proceeding to evaluate the physical realizability of these machines, it is important that I make further precise my notion of ``physical realizability''. The Argument by Infinity and the Argument by Physical Vagueness offered a satisfactory starting point, allowing us to provide arguments against a large number of hypercomputational models. Given that the original statement of the Church-Turing thesis does not \emph{directly} apply to what can be computed by physical systems in general (but is concerned with what can be determined by effective procedures), and given that Copeland and Shagrir aim to discuss machines that compute non-recursive functions in a broad sense of ``compute'' (especially due to the substantial number of hypercomputational models in the literature), the question is: Are there more general constraints on physical computation by which to judge a particular model's physical realizability? \newline
\subsubsection{General constraints on (humanly useful) physical computation}
\indent Piccinini \cite[\S 2]{p-phys} offers a suitable notion of physical computation, which we can use, and can be formulated as the following ``usability constraint'' on physical computation: 
\begin{use-constr}
If a physical process is a computation, it can be used by a finite observer to obtain the desired values of a function.
\end{use-constr}
The statement of the usability constraint can be further clarified if one clarifies what is meant by ``finite observer'' and ``can be used by a finite observer''. A finite observer, in Piccinini's sense, is a human being or any other intelligent being of similarly bounded capacity (he actually proposes another broader definition, but the current definition is sufficient for the purposes I intend). As to the issue of what it means for a physical process to be useable by a finite observer, Piccinini formulates the following six sub-constraints, which shall be summarized here:
\begin{enumerate}
\item \textbf{Readable Inputs and Outputs:} The inputs and outputs of a computation must be readable, in the sense that the inputs and outputs can be measured to the desired degree of approximation. For example, having an infinite precision real number be the input and/or output would not be permitted.
\item \textbf{Process-Independent Rule:} In a genuine computation, the problem being solved (equivalently, the function being computed) must be definable \emph{independent} of the process of solving it (or equivalently, computing it). This sub-constraint entails that modeling one physical system using another is not enough for computing.
\item \textbf{Repeatability:} For a physical process to be a genuine computation, it must be repeatable by any competent finite observer who intends to obtain its results, in the sense that that the same sequence of states can occur either within the same physical system at different times or within two relevantly similar physical systems (e.g. two systems satisfying the same equations).
\item \textbf{Settability:} An ordinary computing system can compute any value of one or more functions, within its limits (e.g. a universal machine can compute the value of any computable function until it runs out of memory or time). Moreover, a system must be \emph{settable}, namely, a user can choose which value of the function is going to be generated in a given
case.
\item \textbf{Physical Constructibility:} If a system cannot be physically constructed, it may count as performing notional computations, which are irrelevant for physical purposes. This fifth sub-constraint is certainly \emph{not} satisfied by a machine that is either affected by the Argument by Infinity or the Argument by Physical Vagueness.
\item \textbf{Reliability:} While the requirement that machines should never break down is unrealistic, it is a requirement that machines be reliable in the sense that they completely operate correctly long enough to yield correct results at least some of the time. A machine that never completes its computation successfully is certainly not worth building.
\end{enumerate}
As noted by Piccinini, the usability constraint is weaker than the verifiability constraint that is discussed and rightfully dismissed by Shagrir and Pitowsky \cite[pp. 90-1]{sp}. Additionally, the notion of finite observer that was defined above makes the usability constraint relevant to \emph{humanly useful} physical computation -- a proper subset of physical computation in general. Our interest is in characterizing humanly useful physical computation because one main goal of the hypercomputation movement has been to argue for the physical realizability of models that compute the non-Turing computable and could, in principle, be used by humans to compute such functions. And the above list of six sub-constraints, though certainly not exhaustive, is enough to exemplify the properties that are required to produce a useful physical computation. Moreover, the use of this usability constraint in \S 4 will further demonstrate its importance. \newline
\indent Admittedly, these contraints are preliminary, more work will have be done in order to make them less open-ended. However, for the extensive counterexamples that Copeland and Shagrir present, it turns out that the Argument by Infinity and the Argument by Physical Vagueness suffice for sharpening the fifth sub-constraint. In other words, if a model violates any of these, then it can be said that it is not ``physically constructible'' (given current physical theories). On the other hand, just because a model is physically constructible does not mean, in my terminology, that it is ``physically realizable''. If the machines that Copeland and Shagrir \cite{copeland-shagrir} present which fall outside the class of Gandy machines do not satisfy this usability constraint (namely by violating at least one of the six aforementioned sub-constraints), then they are ruled out as ``physical realizable'' (I will refer to them as ``physically unrealizable hypercomputers'', or ``physically unrealizable'' for short, or any other applicable synonym). I should note that my criterion encompasses Piccinini's criterion for distinguishing between what he calls ``spurious'' and ``genuine'' hypercomputers, where the former fails to satisfy at least one of the first four sub-constraints on physical computation, and the latter satisfies at least the first four sub-constraints. \newline
\subsubsection{Partially random machines}
\indent Setting aside the case of analog hypercomputation for the moment, the first model that Copeland and Shagrir present is the ``partially random machine'', introduced by Turing \cite{turing1948} in an unpublished (in his lifetime) manuscript entitled ``Intelligent Machinery'' as follows: \newline
\begin{quote}
It is possible to modify the above described types of discrete machines by allowing several alternative operations to be applied at some points, the alternatives to be chosen by a random process. Such a machine will be described as `partially random'. If we wish to say definitely that a machine is not of this kind we will describe it as `determined'. Sometimes a machine may be strictly speaking determined but appear superficially as if it were partially random. This would occur if for instance the digits of the number $\pi$ were used to determine the choices of a partially random machine, where previously a dice thrower or electronic equivalent had been used. These machines are known as apparently partially random. Turing \cite[pg. 416]{turing1948} \newline
\end{quote}
The authors point out that the outputs of a partially random machine may form a non-computable sequence, since as Church \cite{church40} argues, if a sequence of digits $a_1, a_2, \ldots, a_i, \ldots$ is random, then there is no function $g(i) = a_i$ that is computable. However, a random output cannot be used to perform any interesting physical ``hypercomputation'', such as solving the halting problem, since given the usability constraint, there are several reasons to deny that a genuine random process $R$ should be considered a \emph{physical} computation. As Piccinini \cite[pp. 750-1]{p-phys} mentions, there is no process-independent way to define the function $f: \mathbb{N} \to \{0,1\}$ without referring to $R$ itself, thereby violating the second sub-constraint; random processes are not repeatable, thereby violating the third sub-constraint; and lastly, random processes are not settable (which violates the fourth sub-constraint) since the ``user'' of $R$ cannot select the value of $f$ he or she intends to generate. Of course, random processes can be exploited by a computation, but random processes cannot be used to generate the desired values of a function or solve the desired instances of a general problem, such as the halting problem. Physical unrealizability aside, it is not known, as Copeland and Shagrir mention, whether it is possible for there to be partially random (ideal) machines that offer any interesting counterexamples to Thesis M. Even this possibility is unlikely, for in the case of ideal machines where one \emph{can} refer to the possibility of enumeration of infinite sets and the computation of infinite sequences, Theorem 2 of the classic paper by de Leeuw, et al. \cite{probmach} essentially states that anything that can be done by a machine with a random element (say, for the sake of example, it is a binary random device with probability $p$ of producing a 1 in its output) can be done by a deterministic machine (given that $p$ is a computable real number\footnote{In fact, given Turing's description of partially random machines, it seems that the means by which the machine is to make its choices in order to appear it has a random element makes the restriction that $p$ be a computable real number rather reasonable.}). \newline
\subsubsection{Gandy machines that produce a non-computable output}
\indent Of course, the partially random machine is not a Gandy machine; though, the authors also consider machines that do satisfy Gandy's principles but may produce a non-computable output from computable input. It is worth mentioning though that the underlying assumption being made is that this occurs if these machines are embedded in a universe whose physical laws have non-computability built into them, and therefore the authors place themselves in an imagined universe. Hence, almost immediately, this violates the fifth sub-constraint of physical constructibility. \newline
\indent Not to mention that the example they consider are asynchronous networks of Turing machines, where each Turing machine has a timing function $\Delta_k: \mathbb{N} \to \mathbb{N}$ that represents the units of time between its execution of its $i$-th action and its $(i+1)$-th action. A non-recursive source would be needed to be used by the timing function in order for the asynchronous network of Turing machines to compute the non-computable. One way to do this would be to use a real-valued quantity, such as a particle with a non-recursive charge, but the physical drawback that comes with this method is the requirement of arbitrary precision measurement, which is impossible to attain physically. Another approach of accessing a non-recursive information source may be by way of the observation of a non-recursive physical process. Particle decay is examined as a possible approach by Ord, but it is readily pointed out that while it may be physically harnessable, it is ``completely useless for computing definite mathematical results that could not be computed with a Turing machine'' \cite[pg. 26]{ord}. Thus, the question arises as to whether there exist other non-recursive information sources that may be \emph{mathematically} harnessable and useful. One possibility is if the bits of a constant such as $\pi$ or $\tau$ were given by some natural process, but there is little possibility that such a natural process exists due to the Argument by Infinity, thereby violating the fifth sub-constraint in another way. \newline
\subsubsection{Trial and error hypercomputation}
One notion of discrete hypercomputation not covered by Copeland and Shagrir is hypercomputation by means of trial and error. Kelly \cite{kelly} writes that Turing machines (with finitely many states, rules, and symbols) can effectively compute non-recursive sets if they are allowed retractions. Thus, the computation produces output from time to time and is guaranteed to eventually produce the desired correct output, but with no bound on the time needed for this to happen. Once the correct output has been produced, any subsequent output will just repeat this correct result. This idea is certainly not new, originally pioneered by Gold \cite{gold} and Putnam \cite{putnam}. But as Davis points out, there is the obvious setback that ``someone who wishes to know the correct answer would have no way to know at any given time whether the latest output is the correct output'' \cite[pg. 128]{davis-sr}. Hence, trial and error hypercomputability falls to the Argument by Physical Vagueness, therefore violating the fifth sub-constraint. From this discussion, we naturally segue to supertasks. \newline
\subsubsection{Supertasks as hypercomputation}
\indent The final set of machines that Copeland and Shagrir discuss are supertasks, in particular the accelerated universal Turing machine (AUTM) and the approach of using supertasks in relativistic space-time. In the early 20th century, Blake \cite{blake}, Weyl \cite{weyl}, and Russell \cite{Russell} independently observed that a process that performs its first step in one unit of time and each subsequent step in half the time as the preceding step, can complete an infinity of steps in two units of time. A universal Turing machine working in this manner is called an accelerated Turing machine. This device is capable of deciding the halting problem; not to mention that the definition of the Turing machine does not specify how long an individual step may take during computation, so the argument is that the accelerated Turing machine does not conflict mathematically with the classical notion (cf. Calude and P\u{a}un \cite{c-p}). \newline
\indent The AUTM is an example of a machine obeying Newtonian mechanics \cite{copeland-acc} and hence falls outside of the class of Gandy machines since it does not obey Principle IV. Specifically, it violates Gandy's requirement that ``that there is an upper bound (the velocity of light) on the speed of propagation of changes'' \cite[pg. 126]{gandy}. The authors also include Davies' \cite{davies2001} proposal for the construction of an infinitely accelerating Turing machine which involves an infinite series of machines that, after performing some computational operations, builds a smaller and faster copy of itself, and so on, infinitely so. Again this proposal is consistent with Newtonian mechanics, for Davies' goal was to produce machines similar to Babbage's Analytical Engine, ``consistent with any physics known in the year 1850'' \cite[pg. 672]{davies2001}. Although this does not break the upper-bound on the velocity of light, it still does not obey Principle IV as it depends on an absence of a lower bound on the size of the atoms. Not only is this latter proposal evidently susceptible to the Argument by Infinity, but a relevant question is why would one want to conform to what was known about physics around 1850? Precisely where infinities arise in Newtonian frameworks for hypercomputation and where Newtonian mechanics breaks down are at extremes such as the velocity of light or on the atomic scale. Even as Davies admits, the quantum mechanical nature of matter makes his construction physically impossible. \newline
\indent Thus, I shift my focus to hypercomputation by way of modern physical theories, and inquire as to whether or not this is physically realizable. The supertasks in relativistic space-time that Copeland and Shagrir present are unlikely to be physically realizable (cf. Piccinini \cite[\S 4.2]{p-phys}). It is worth noting that Piccinini considers these PMH (Pitowsky-Malament-Hogarth) relativistic hypercomputers to be ``genuine'' hypercomputers, as they satisfy sub-constraints 1-4, but at the cost of their physical constructibility (sub-constraint 5) and reliability (sub-constraint 6) being unlikely. \newline
\indent But one can say more general things about supertasks as hypercomputation. Shagrir \cite{shagrir-acc04} asserts that simply because the AUTM performs a supertask, does not mean it can compute what Turing machines cannot compute; instead, what matters is the \emph{computational structure} involved. Since the AUTM and the Turing machine have the same structure, their computational power is equivalent. As Shagrir \cite{shagrir-acc04} argues, the AUTM does not even solve the halting problem although it performs a supertask. However, this does not mean that supertask machines cannot be hypercomputers (but their hypercomputational capabilities have less to do with their performance of a supertask and primarily has to do with their structure differing from that of the Turing machine). As examples of this, in a later paper, Shagrir \cite{shagrir11} considers as supertask hypercomputers the AUTM$^{+}$ and relativistic hypercomputation. Regarding the latter, I have already mentioned that the PMH proposal is not physically realizable, due to its not satisfying the usability constraint. Regarding the AUTM$^{+}$, the difference between the AUTM and the AUTM$^{+}$ has to do with the inclusion of a second moment limit stage in the latter such that one of the machine's stages is not completely determined by the configuration of the previous stage (unlike a Turing machine). The reason, as Shagrir clarifies, is that ``there are infinitely many stages between any given stage that precedes the second moment and the second moment limit stage'' \cite[pg. 54]{shagrir11}. Hence, this not only violates Gandy's Principle I, but is also not physically realizable, as it violates the fifth sub-constraint due to the Argument by Infinity.\newline
\subsubsection{Comparative notions of deterministic machines}
\indent Moreover, Copeland and Shagrir \cite{copeland-shagrir} argue that PMH is of interest because it \emph{prima facie} complies with Gandy's Principles I-IV, yet computes the non-computable. Namely, PMH is a discrete state machine that consists of two communicating digital computers and PMH is consistent with the laws of general relativity, yet there are two important differences which Copeland and Shagrir are keen to point out. First, Gandy assumes that processes that consist of infinitely many steps do not terminate, whereas PMH allows terminating process that consist of infinitely many steps. Second, if one adheres to Gandy's preliminary definition of deterministic that ``the subsequent behaviour of the device is uniquely determined once a complete description of its initial state is given'' \cite[pg. 126]{gandy}, then PMH is deterministic in this sense. But, the more mathematically precise definition of what Gandy means by ``deterministic'' presented in Principle I explicitly requires that the configuration of each state is uniquely determined by the configuration of the previous step; PMH does not follow Principle I, and thus is not deterministic in Gandy's sense. Copeland and Shagrir argue that PMH can be considered deterministic in a sense that is not ``Gandy-deterministic'' based on the initial state of the machine, and that ``this sense of determinism is in good accord with the physical usage whereby a system or machine is said to be deterministic if it obeys laws that invoke no random or stochastic elements'' \cite[pg. 228]{copeland-shagrir}. However, if the PMH proposal is unlikely to be physically realizable (since it does not comply with the usability constraint), then Gandy's definition of determinism cannot be justifiably ruled as too narrow on the grounds of not being ``in good accord with the physical usage.'' \newline
\subsubsection{Analog hypercomputation}
\indent Perhaps what is more interesting is the physical possibility of analog hypercomputation (which Copeland and Shagrir also address), since as mentioned at the beginning of \S 2, the discussion of physical properties (which can be used to possibly construct an analog hypercomputer) that cannot be simulated by Turing machine has spanned over at least four decades, even before the ``hypercomputation movement'' began. I have hitherto discussed in this section the impossibility of discrete deterministic hypercomputers existing in the real world based on Gandy's principles, and have ignored the issue of analog hypercomputers for the time being, primarily due to the fact that Gandy did not discuss the analog case in his original paper \cite{gandy}. \newline
\indent Such proposals for instance that I have excluded from consideration include those by Kreisel \cite{kreisel1974}\cite{kreisel1982}, Pour-El and Richards \cite{per79}\cite{per89}, and da Costa and Doria \cite{cd91}. Due to a lack of space, I will not go into the details of these proposals, especially since several of them (in particular the Pour-El and Richards results) have been well-known for quite some time. Instead, my purpose is to formulate a response to the physical realizability of analog hypercomputation in general. In fact, in an unpublished hand-written manuscript, Gandy \cite{gandy1993} approaches analog machines (which were excluded in Thesis M) and discusses examples of physical systems (including the aforementioned proposals), both classical and quantum mechanical, for devices that compute the non-computable. His argument is that given a reasonable definition of ``analog machine'', such a machine cannot calculate the values of a number-theoretic non-computable function, upon accepting a continuously variable computable input. Gandy's position is summarized as follows (cf. Kieu \cite[pp. 557-8]{kieu-qhyp} and Copeland and Shagrir \cite[pp. 221-2]{copeland-shagrir}): \newline
\indent If we let $A$ be a recursively enumerable non-recursive set (for example the set which represents the halting problem), then a decision problem can be posed as the question of whether or not an integer $j \in A$. Because it is recursively enumerable, there exists a total computable function $a: \mathbb{N} \to \mathbb{N}$ which enumerates $A$ without repetitions. Next, define the \emph{waiting-time} function $\nu$ as
\begin{equation*}
\nu(j) := \mu n[a(n) = j],
\end{equation*}
which is a partial recursive function whose domain is $A$, and which is not bounded by any total computable function. Moreover, $\nu(j)$ gives the least $n$ (``time'') upon which $j$ is confirmed to be a member of $A$. Now, Gandy argues that for any particular analog machine there is an upper bound $J$ on the inputs it can accept, and he then proceeds to define the total function $\beta(J)$ as
\begin{equation*}
\beta(J) := \operatorname{max}\{\nu(j): j < J \mbox{           } \& \mbox{           } j \in A\}
\end{equation*}
(with $\operatorname{max} \emptyset = 0$). $\beta(J)$ is not computable, and indeed it eventually majorizes every computable function. Gandy's claim is that with given $J$, one cannot design an analog machine (whose behavior is governed by standard physical laws) which will give correct answers to all the questions $j \in A$? for $j < J$ unless one knows a bound $B$ for $\beta(J)$. To illustrate the significance of the wording of the claim, suppose (what is quite plausible) that someone proves that $j \not\in A$ for all $j < J = 10$; then he can design a machine which always outputs ``NO'' for $j < J$. But, because of this proof he does in fact know that $\beta(J) = 0$. Of course, if one knows a bound $B$ as above, then one does not need an analog machine to settle the decision problem. One simply computes $a(n)$, for all $n < B$, to see if $j$ is obtained as a value, and thus belongs to $A$, or not. \newline
\indent It is worth mentioning that Gandy's argument, though rather encompassing, was intended as a ```challenge rather than a dogmatic assertion''' \cite[pg. 221]{copeland-shagrir}. Kieu attempts to take on this challenge by arguing that Gandy's conception of an analog machine is too narrow, stating \newline
\begin{quote}
Gandy also recognised that his claim might not stand up to the kind of analogue machines whose behaviour depends on a single quantum (similar to our process involving Hamiltonians which have discrete spectra of integer values). Neither would the claim be valid if the physical theory involved has elements of non-computability already built-in. But here we have argued all along that Quantum Mechanics with its measurement postulate has, in a sense, some elements of non-computability through its \emph{intrinsic randomness}. To wit, randomness is `noncomputable' as it cannot be algorithmically computed....Furthermore, our proposal is not subject to Gandy's treatment since ours is a kind of probabilistic computation, where the probability that we get a wrong answer can be made arbitrarily small but is not exactly zero. Kieu \cite[pg. 558]{kieu-qhyp} \newline
\end{quote}
It is interesting further to mention that Copeland and Shagrir mention Kieu's counterarguments to Gandy's analog machine claim, writing that ``we will not pursue these matters here, but observe that the deep question of whether there are (ideal) analogue machines which satisfy usability constraints and which, when provided with a computable input, will give a non-computable output,
is open'' \cite[pg. 222]{copeland-shagrir}. However, certainly Kieu's proposed hypercomputer \cite{kieu2001} \cite{kieu2003} does not satisfy the usability constraint we are using, for it certainly does not satisfy the fifth sub-constraint, due to the Argument by Infinity (see, for example, Hagar and Korolev \cite{h-k}). And if one is to accept Kieu's proposal based on occupation numbers (which can be at least ${{10^{10}}^{10}}^{10}$ (cf. \cite[pg. 5]{davis2006}), so when written out, their decimal expansion is larger than the diameter of the galaxy) then another sub-constraint, namely the first sub-constraint, is violated since the output may not be readable. For now, it is a very strong possibility that quantum mechanics (the weakness which Kieu alleges Gandy's claim to be susceptible to) is unlikely to allow for physically realizable hypercomputation (even despite the fact that, Kieu's proposal aside, quantum operations call for the execution of primitive operations that Turing machines cannot perform); instead, a more interesting prospect is that it may provide great future aid in solving problems of current practical intractability, which I will discuss in \S 4.2.1 (as well as the more general evidence regarding the \emph{limits} of quantum computation). \newline
\subsubsection{Gandy's reductio ad absurdum}
\indent In sum, as can be seen, unlike with Turing's and Gandy's definitions of computability, for which one can obtain a precise axiomatization, there is currently no single set of principles under which all hypercomputational models (provably) fall. Of course, even if there is no single set of principles, all it would take is one physically realizable notion of hypercomputation that goes beyond the universal Turing machine to justify this idea; nonetheless, given the usability constraint, the arguments in this section have supported Gandy's principles as sufficiently general for real machines. Instead, one attains some degree of mathematical precision in terms of understanding the question ``What is hypercomputation?'' as a weakening or contravention of at least one of Gandy's four principles (at least for non-analog hypercomputers, the analog case was dealt with separately). \newline
\indent Regarding the issue that Copeland and Shagrir take with Gandy's \emph{reductio ad absurdum} argument, namely, as Israel puts it, that each of Gandy's principles is necessary to avoid the ``absurdity or vacuity...[that] every number-theoretic function is computable'' \cite[pg. 197]{israel}, they assert that each number-theoretic function \emph{is} computable, relative to some set of capacities and resources. As a result, they do not agree that there is any absurdity in allowing that every number-theoretic function is computable (relative to some index), especially since questions about which functions are computable relative to certain physical theories are seldom trivial. ``When classicists say that some functions are absolutely uncomputable, what they mean is that some functions are not computable relative to the capacities and resources of a standard Turing machine. That particular index is of paramount interest when the topic is computation by effective procedures. In the wider study of what is computable, other indices are of importance'' \cite[pg. 229]{copeland-shagrir}. \newline
\indent But, my arguments did not appeal to the capacities and resources of a standard Turing machine since the \emph{gedankenexperiment} of the Turing machine was not originally intended for physical computation but regarded the limits of effective procedures. Instead, my assertions were based on the epistemic usability constraint which is reasonable given current, modern physical theories. Thus, Gandy's axiomatic principles for mechanisms appear to be general enough to make hypercomputation physically unrealizable, and for now, given the arguments presented in this section, it appears that what Piccinini calls the Modest Physical CT\footnote{``Modest'' with respect to the Bold Physical CT, which states that \emph{any} physical process is (Turing-)computable.} holds, namely, that any function that is physically computable is (Turing-)computable, for which Piccinini already provides extensive support. \newline
\indent As of now, the physical theories we have deem hypercomputation physically implausible. What may be required to build a hypercomputer is a groundbreaking new physical theory. Now, it is important to mention that I am not saying that the theoretical pursuit of hypercomputation and the models presented so far are to be discarded, since some may have interesting theoretical applications (e.g. the $o$-machine which was seminal to the notion of uniform relative computability as well as the PMH proposal which is an intriguing contribution to debates on the foundations of physics); however, in the case of physical realizability, these models are not to be considered for that purpose. Therefore, to make hypercomputation relevant to the problems and applications of today, a shift in emphasis must occur from computing the non-computable to expressing/modeling certain processes in a superior manner than the closed-system paradigm of the Turing machine. The justification for this shift will be explained in the next section.

\section{The practical limitations of Turing-equivalent machines: A proposal}
\subsection{Computational theory vs. Computational practice}
\indent While the hypercomputation movement has claimed that the theoretical limits of computing are physically realizable, complexity issues in computational practice have not been given much focus (if at all) by this movement. As an immediate important example, which Davis \cite{davis-sr} brings to light, the work of Cook and Levin on NP-completeness made severely relevant the issue of problems for which feasible algorithms are not known and for which Turing computability does not help, since the number of steps required to solve the problem grows exponentially with the length of the input, making their use in practice a problem. Certainly the lofty, to say the least, attempt at computing the non-computable in practice is far from the already difficult problem of solving intractable problems in practice. \newline
\indent As Feferman \cite{feferman-o} distinguishes: The aim of computational practice is to produce hardware and software that is, among other things, reliable, efficient, and user-friendly. The aim of the theory of computation, through its employment of logic and mathematics, is to provide limits to what is feasible and to provide a body of concepts around which to organize experience and a body of results predicting correctness, efficiency, and versatility in order to aid engineers in producing hardware and software that meet the aforementioned requirements. Certainly this comparison coincides with the usability constraint, which I used in \S 3 to evaluate the possibility of physical hypercomputation more generally. While the focus of this paper so far has either implicitly or explicitly been on the universal Turing machine, its use in computational practice is questionable, and Turing-equivalent models are often more applicable to a certain area than Turing's $a$-machine, particularly because the Turing machine was not intended to be a practical computing device. Nonetheless, there still remain aspects of computational practice for which Turing computability does not help, and this will underlie my proposal for a more plausible hypercomputability. And in order to explicate and justify my assertions when examining the practical limitations of Turing-equivalent machines, I will appeal to the usability constraint. \newline
\indent  Many machines that are thought to have a smaller instruction set, use less memory, or are more applicable to a certain area than a Turing machine can be shown to have at most the same computational capacity as the universal Turing machine \cite{hopcroft-ullman}. It is important to note that Turing's $a$-machine was not meant to model a ``real'' computer, but to determine the limits of effective procedures, and of course, a computer with fixed internal storage was developed only \emph{after} Turing's influential paper. In fact, it is well-accepted that the Turing-machine itself does not provide a realistic model of actual computers and that register machines more closely resemble modern digital computers \cite[\S 3.5]{soare}. \newline
\indent I am in agreement with Feferman that the notions of relative computability have a much greater significance in the realm of computational practice than those of absolute computability. It is here that Turing's $o$-machine (which certainly is not a ``forgotten idea'' of Alan Turing, as Copeland and Proudfoot \cite{copeland-proudfoot} claim) becomes of prime relevance. Feferman attributes this to the modular organization of hardware and software, citing the example of built-in functions in computers as ``black boxes'' (or ``oracles'') in practice. He also mentions in a footnote that several people have suggested to him interactive computation as an example of Turing's oracle in practice, but ``while I agree that the comparison is apt, I don't see how to state the relationship more precisely'' \cite[pg. 340, fn. 8]{feferman-o}. Given that this was written in 1992, and since then technology has greatly expanded and is now significantly centered around the World Wide Web, online computation is the \emph{most} pervasive example of Turing's oracle in day-to-day computational use. And, in fact, in a 2009 article, Soare \cite{soare09} does provide a precise characterization of the relationship between Turing's oracle and online computing. It is important to note that the Turing machine itself does not communicate with its environment (i.e. whatever is external to it) during computation. But, on the other hand, online processes do involve interaction with the environment, not to mention that the relatively commonplace notion of upgrading hardware or software on a computer offers new functionalities, and as a result, its programs change due to external influence. This is problematic for the Turing machine paradigm in terms of encapsulating such processes; therefore, the notions of absolute computability do not capture certain, rather ubiquitous, aspects of modern computing. None of this requires one to go significantly beyond the theoretical framework of what was developed by the ``classicists'', since \newline
\begin{quote}
The theory of relative computability developed by Turing and Post and the $o$-machines provide a precise mathematical framework for database or online computing just as Turing $a$-machines provide one for offline computing processes such as batch processing. Oracle computing processes are those most often used in theoretical research and also in real world computing where a laptop computer may communicate with a database like the World Wide Web. Often the local computer is called the ``client'' and the remote device the ``server''. Soare \cite[pp. 370-1]{soare09} \newline
\end{quote}
When it comes to computational practice, however, there are still important developments yet to be made regarding interactive computation that are not accounted for by the theoretical models, that I will discuss later on in this section. \newline
\subsection{Towards a plausible extension of the classical paradigm}
\indent Just as one refers to the ``theory and practice of computability'', I will similarly refer to the ``theory and practice of hypercomputability''. With regard to the \emph{theory} of hypercomputability, I do not propose any drastic changes. I do not think that models of hypercomputation that supposedly compute the non-computable should be discarded since they may be of theoretical interest (not to mention that the theoretical study of non-computability has been an active field of research since the 1940s starting with Kleene and Post). On the other hand, concerning the \emph{practice} of hypercomputability, it is the purpose of this final section to justify a direction different from the hypercomputation movement's current practical aims. Instead of attempting to provide some physical credence to computing the non-computable, a reasonable extension to the classical paradigm should adequately encapsulate certain modern computational paradigms for which Turing computability does not help, such as those that are computationally intractable in practice, as well as interactive and online computing. \newline
\indent Before specifying my proposal, it is relevant to mention that Gurevich \cite[pg. 8]{gurevich} asserts that most real world machines do not satisfy Gandy's principles, but he provides very little justification for his claim, stating only that ``the burden of proof is on the proponents of the approach.'' While my arguments in \S 3 provided extensive support for the claim that Gandy's principles for mechanisms are not too narrow, my proposal's acceptance should not be entirely contingent on my arguments in \S 3. They do not settle, for instance, the question of what a (not necessarily humanly useful) physical computation is, and one could argue that a clear definition of this notion is needed in order to satisfactorily argue against hypercomputation as being physically realizable in a general setting. Surely, one can raise other philosophical questions (which may or may not be related to the physical intepretation of computation) and use that to possibly point to gaps in the argumentation. Instead, \S 3 should be viewed as a \emph{first step} towards challenging hypercomputability from a more general perspective, an improvement over the weaker model-specific arguments presented thus far. And certainly these more general arguments should be investigated further. All this said, there is also a trivial, but crucially important observation to make: nobody has built a hypercomputer or is even close to building one! Thus, a more plausible proposal for extending the classical paradigm is needed, given that hypercomputation has not been of much use to computer science practice (which is evident even from the above mentioned trivial observation). My proposal, in essence, is to understand and investigate modern-day computations that are not adequately encapsulated by the closed-system paradigm of the Turing machine. It dispenses with the notion that the classicists (G\"{o}del, Church, Turing, Post, etc) figured it ``all'' out, and is instead meant to embrace the more reasonable observation that computer science is ever-evolving and that even our notion of computation is not crystallized (a major theme, for instance, of Gurevich's paper \cite{gurevich}). \newline
\indent I will proceed to specify my proposal below. This specification will inevitably entail a mix of survey and speculation, in contrast to the detailed argumentation provided in \S 3. \newline
\subsubsection{Intractability in practice}
\indent Even if a hypercomputer is built, there is still the question of the time it takes to provide answers. If there is significant delay, such as if the time it takes to produce an output is exponential with respect to the size of the input, then reliability (sixth sub-constraint) becomes an issue if the inputs it is given (for some particular purpose) are large enough that the output is never returned in a finite amount of time that is useful for practical purposes (such as 1000 years, $10^{10}$ years, etc.). However, one does not need hypercomputation to realize that complexity is a relevant issue in the real world -- these issues are a well-known occurrence in current computational practice. The question that arises is, how can one effectively solve certain important problems that are intractable? A major step forward would be to change the computational paradigm used in practice from the commonplace electronics (which has been a major goal of the natural computing movement). \newline
\indent Quantum computing, especially, may help to overcome calculations that are computationally infeasible (not to mention, as alluded to in \S 3.2.7, that not all of the primitive operations of a quantum Turing machine, introduced by Deutsch \cite{deutsch}, can be performed by a person simply using pencil and paper -- but it should be noted that it cannot compute non-recursive functions). Already this is underway, the most famous example being Shor's algorithm \cite{shor} in 1994, which demonstrated that the factoring of large integers could occur in polynomial time on a quantum computer, a task insurmountable for the standard electronic computer (and particulary relevant to cryptography, for example). Quantum computers that are able to do this have yet to be built; hence, the fifth sub-constraint of physical constructibility is not yet satisfied, though their realizability is certainly not impossible (an excellent discussion is provided in Nielsen and Chuang \cite[Ch. 7]{quantum-textbook}). The field of quantum computing has many opportunities for new innovations, and thus one goal of a feasible extension to the classical paradigm could be to make such quantum computations physically realizable. Despite this initial optimism, there is also evidence regarding the limits of quantum computing (cf. Nielsen and Chuang \cite[Ch. 6.7]{quantum-textbook}). This is based on the following black-box abstraction of what it means to a solve a problem in NP. There is a boolean function $f(.)$ that takes an $n$-bit string as an input. We can evaluate $f(.)$ any input of our choice, and we are to determine if $\forall x. \mbox{              } f(x) = 0$ or if $\exists x. \mbox{       } f(x) = 1$. The end result is that in computing boolean functions using a black box, quantum algorithms may provide only a polynomial speedup over classical algorithms \emph{at best} and even this is not possible in general. No exponential speedup over classical algorithms is possible without additional information about the \emph{structure} of the black box, which means that this does not rule out, for instance, a quantum polynomial time algorithm for an NP-complete problem such as 3-SAT. However, this cannot come just from understanding the quantum model, but from an algorithmic breakthrough on the specific problem itself.  My mention of quantum computing is intentionally brief as its potential is rather well publicized, and I will discuss some other lesser known alternatives to electronics that show interesting potential in the realm of tractability. \newline
\indent Introduced by P\u{a}un \cite{paun2000}, P systems abstract distributed parallel computing models from the structure and function of a living cell. There are many variants of P systems, but the basic model of P systems is a membrane structure consisting of several cell membranes which are hierarchically embedded in a main membrane, termed the skin membrane. Membranes delimit regions, wherein objects are placed. Objects evolve according to evolution rules. These rules are applied to objects in the region with which they are associated and can modify the objects, send them outside the current membrane (or to an inner membrane), and can dissolve the membrane. The whole functioning of the system is governed by evolution rules which are applied in a maximally parallel manner, namely, at each step, all objects that can evolve should evolve. A \emph{computation} starts from an initial configuration and proceeds using the evolution rules, and the computation is considered completed when it ``halts'' -- no further rule can be applied to the object present in the last configuration. P systems are both capable of Turing universal computations and, in cases where enhanced parallelism is provided, can solve intractable problems in polynomial time. As an example of the latter, P\u{a}un \cite{paun-np} proposes a class of P systems where evolution rules are associated with both objects and membranes, and membranes can not only be dissolved, but can multiply themselves by division. In this way, the number of membranes can grown exponentially, and by making use of this enhanced parallelism, P\u{a}un proves that one can compute faster, namely, the well known NP-complete SAT problem can be solved in polynomial (linear) time (more on the computational efficiency of P systems can be found in a survey article of P\u{a}un \cite[\S 8]{paun-surv}). The computational universality of P systems and the possibility of solving NP-complete problems in linear time makes membrane computing a promising candidate as a model of computation. But it is important to mention that there it is still an open question of implementing a P system of this type in either biochemical or electronic media; therefore, there is still work to be done regarding whether or not this P system can fully satisfy the fifth sub-constraint of physical constructibility. An informative article briefly overviewing the fast occurring developments in P systems is provided by P\u{a}un \cite{paun-surv}, and in \S\S 9 and 10 of this article, he discusses current research topics in the study of P systems as well as their applications. \newline
\indent Computational mathematics may be feasible with DNA (beyond just primitive mathematical operations). The recent field of DNA computing is the study of the information-processing capabilities of DNA, based on the massive parallelism of DNA strands and Watson-Crick complementarity. The key advantage to DNA is that it can make computers much smaller than before, while at the same time maintaining the capacity to store prodigious amounts of data. For instance, computationally intractable problems can be solved by an exhaustive search through all possible solutions. But the obvious problem with this approach for modern computers is that the search space is too large. On the other hand, the density of information stored in DNA strands and the ease of producing many copies of them might allow these brute force searches to be possible in practice. Since Adleman's \cite{adleman} pioneering paper, DNA computing has become a rapidly evolving field with its primary focus on developing DNA algorithms for NP-complete problems. Unlike quantum computing in recent years, though, the viability of computational mathematics on a DNA computer has not yet been fully demonstrated. In fact, only recently have the primitive operations in mathematics (i.e. addition, subtraction, multiplication, and division) been fully implemented and optimized. A recent paper of mine \cite{nayebi} discusses the theoretical possibility of implementing fast matrix multiplication algorithms (Strassen's algorithm is the focus\footnote{There have been several recent improvements to the exponent on the theoretical runtime of matrix multiplication (by Stothers \cite{stothers} and Williams \cite{williams}), which is conjectured to optimally be $O(n^2)$ for the multiplication of two $n \times n$ matrices. But unlike Strassen's algorithm \cite{strassen}, these algorithms are not used for matrices of practical sizes ($n < {10}^{20}$, \cite{nayebi}).}) with DNA due to its massive parallelism, and fast matrix multiplication serves as a bottleneck for many algorithms (such as inversion, computing determinants, and graph theory). \newline
\indent Of course, there is no current advantage that DNA has over the electronics with regards to mathematical computing, not to mention the question of determining reliability \cite{losseva}, and at the moment DNA computing has merged to controlling and mediating information processing for nano structures and molecular movements. It is also relevant to note that Sieg \cite{sieg2000} mentions that the Gandy machine is an abstract mathematical definition that embodies generally accepted properties of parallel computations (and indeed, with John Byrnes, he formulates an abstract model for parallel computations in \cite{sb1999b}), and that they may be simplified via a graph-theoretic or a category-theoretic representation. However, ``what is needed most...is their further mathematical investigation, e.g., for issues of complexity and speed-up, and their use in significant applications, e.g., for the analysis for DNA computations or of parallel distributed processes'' \cite[pg. 402]{sieg2000}. Whether or not Gandy machines are the most ``favorable'' model of parallel computations is a separate matter, but Sieg's consideration does bring to light the question of complexity issues in DNA computation. Hence, another goal of a reasonable extension to the classical paradigm could be to address the issues concerning reliability with DNA computations (which is relevant to the sixth sub-constraint of the usability constraint), as well as understanding with mathematical rigor an investigation of complexity issues that arise in DNA computations. And an overall long-term goal would be to pursue other computational paradigms that would allow for intractable problems in practice to potentially become tractable. As a starting point towards this long-term aim, Eberbach and Wegner \cite[\S 5]{eberbach-wegner} consider four guiding principles with which to deal with intractability (and I include the paradigms considered so far here in this list), namely:
\begin{enumerate}
\item \textbf{Interaction:} through feedback with the external world, where the ``world'' can be represented by an environment or other computational agents, and ``feedback'' can be in the form of direct advice or by some performance measure obtained by interaction with the environment describing how good a solution is and/or the resources used to find it. Performance measures are used by evolutionary algorithms, anytime algorithms, and the \$-calculus (``cost calculus''), for example.
\item \textbf{Evolution:} transforming a problem to a less complex or incomplete one. This approach is used by structural programming, approximation algorithms, anytime algorithms, P systems, and the \$-calculus, to name a few.
\item \textbf{Guessing:} selecting randomly a path to find a solution. For instance, this approach is used by nondeterministic, probabilistic, ergodic, and evolutionary algorithms.
\item \textbf{Parallelism:} exploring all possible solutions simultaneously usually by trading time complexity for space complexity. This approach is employed by quantum computing, DNA computing, P systems, and the \$-calculus.
\end{enumerate}
Based on these principles, they present three other alternatives than the ones mentioned here for dealing with intractability: evolutionary computation, anytime algorithms, and the \$-calculus. The reader is referred to the relevant sections (\S\S6.3-6.5) of \cite{eberbach-wegner} for a detailed discussion of these alternatives. Hence, not only are there six broad proposals for making hard computational problems tractable in practice, with each proposal offering many opportunities for further research, but the principles of interaction, evolution, guessing, and parallelism may also be used to serve as the basis of further inquiry into seeking out new methods of successfully dealing with intractability.\newline
\subsubsection{Paradigms of interactive computation}
\indent The above mention of parallel computations and interaction allows for a natural segue into discussing the role an extension to the classical paradigm could play in interactive computing. Before doing that, a few misconceptions must first be avoided. Eberbach, Goldin, and Wegner \cite{egw} define the following:
\begin{superT}
All computation, including that which cannot by carried out by a Turing machine.
\end{superT}
And they further mention three principles that allow them to derive ``super-Turing'' models of computation:
\begin{enumerate}
\item Interaction with the world;
\item Infinity of resources;
\item Evolution of the system.
\end{enumerate}
One specific example of what they call ``super-Turing computation'' is the dynamic interaction of clients and servers on the Internet because the ``resulting concurrent interaction of multiple clients and servers, with a dynamic configuration of communication links and nodes, cannot be described as a Turing Machine computation, which must be sequential, static, and with all input predefined'' \cite[pg. 178]{egw}. Their use of ``super-Turing computation'' is misleading, as it could be misconstrued as allowing for hypercomputation and that given their example of clients and servers, this every-day interaction (in the modern world) could provide a means of physically implementable hypercomputation. But there is a difference between \emph{expressing/encapsulating} processes in a manner superior to that of the Turing machine, and computing more than a Turing machine. While certainly (as I have previously established) the ``closed-system'' paradigm of the \emph{a}-machine does not allow for the encapsulation of interactive processes, it must be made clear that the three principles of ``interaction with the world'', ``infinity of resources'', and ``evolution of the system'' do not allow for physically implementable hypercomputation. In fact, in \S 6.1 of their paper, they present several ways of solving the halting problem via the principles of interaction, infinity, or evolution. However, they remark that ``none of these solutions of the halting problem ($H$) are algorithmic. In the proofs, either some steps of the solutions do not have well defined (and implementable) meaning -- this is when we refer to the help of oracles -- or we use infinite resources (an infinite number of steps requiring perhaps an infinite initialization time)'' \cite[pg. 189]{egw}. Thus, it is clear that none of these ``hypercomputations'' are physically realizable, either due to the Argument by Infinity or the Argument by Physical Vagueness. \newline
\indent It is further interesting that they mention Turing's work relating to the $o$-machine. In fact, a subsequent inquiry into this aspect of Turing's research does reveal a framework for online computing long before the advent of the World Wide Web, therefore making Goldin and Wegner's bold statement rather unnecessary: \newline
\begin{quote}
Turing machines had been accepted as a principal paradigm of complete computation, and it was premature to openly challenge this view in the late 1970s and the early 1980s....A paradigm shift is necessary in our notion of computational problem solving, so it can provide a complete model for the services of today's computing systems and software agents. Goldin and Wegner \cite[pg. 102]{wg} \newline
\end{quote}
What Wegner, Goldin, and Eberbach fail to clarify is the distinction between absolute and relative computability. Clearly, the notions of absolute computability (such as the Turing machine and the Church-Turing thesis) do not properly encapsulate interactive processes. None of this requires a paradigm shift from the theoretical groundwork laid by the ``classicists'', for the mathematically precise notions of relative computability become relevant here. Soare \cite[pg. 382]{soare09} combines the work of Turing \cite[\S 4]{turing1939} and Post \cite[\S 11]{post1944} on relative computability and restates it in modern terms as follows
\begin{pt}
A set $B$ is effectively reducible to another set $A$ iff $B$ is Turing reducible to $A$ by a Turing oracle machine $(B \le_{T} A)$.
\end{pt}
Although technically speaking, the distinction between absolute and relative computability does not really require a separate thesis, Soare's presentation does explicate this distinction rather nicely. He notes that while Turing's brief introduction of oracles in his 1939 doctoral dissertation \cite{turing1939} did not state this as a formal thesis, it is largely implied by his presentation. Furthermore, Post \cite{post1944} in 1944 makes this thesis explicit, claiming that it is the formal equivalent of the notion of \emph{effectively reducible}, ``a step as significant as the Church-Turing characterization of `calculable' '' \cite[pg. 382]{soare09}. Soare further asserts that under the Post-Turing thesis, real world online or interactive processes can be described by Turing's $o$-machine, and that ``discussing \emph{only} Turing $a$-machines in modern texts, or \emph{only} the Church-Turing Thesis, and not the Post-Turing Thesis on oracle computers, is like discussing only batch processing machines of the 1950's long after the emergence of online computing'' \cite[pg. 388]{soare09}. Hence, though more needs to be mentioned in the computer science literature about relative computability (since an asymmetrical emphasis has been placed on the results of absolute computability), it already offers a rigorous mathematical framework for interaction, unlike Wegner's interaction machines \cite{wegner-interact}. \newline
\indent However, the important take-away from all this is that nevertheless, the notion of Turing computability does not suffice for adequately encapsulating interaction. As an example to explicitly illustrate this in relation to modern computing systems, I will consider the Actor model, a concurrent model of computation first introduced in 1973 by Hewitt, Bishop, and Steiger \cite{actor73} (subsequently developed by Clinger \cite{clinger} and Agha \cite{agha}, among others), one of a number of message-passing models developed in the 1970s, along with Hoare's Communicating Sequential Processes (CSP) and Milner's $\pi$-calculus. Conventional models of sequential computation (such as the Turing machine, the $\lambda$-calculus, etc.) employ the notions of \emph{global time} and \emph{global state}, wherein computation is carried out from one well-defined state to the next, and the transitions between global states are linearly ordered in the global time of the system. This global state approach was continued in automata theory for finite state machines and push down stack machines (as well as for their nondeterministic versions). When computation is not sequential, these notions may not be appropriate. Instead, certain concurrent systems (such as multiprocessor systems and geographically distributed systems \cite{clinger}) are better analyzed by splitting the global state into pieces and viewing the net computation as a set of local computations that interact through message passing. The Actor model emphasizes these notions of local time and local state, and in contrast to the boundedness and locality constraints for Turing computors, the Actor model does not carry out computation in one place, from one well-defined state to the next; computation is distributed in space, and ``actors'' are treated as universal primitives of concurrent (digital) computation. Local times are represented by the arrival orderings of actors, which operate independently of each other except when they interact through message passing, and the communications between actors are represented by the activation ordering. From the actor point of view, sequential computations are distinguishable from concurrent computations by their event diagrams, and are in fact a special case of concurrent computations in that the activation ordering of a sequential computation is linear and its event diagram is a conventional execution sequence. \newline
\indent The Actor model is also of use in computational practice -- an example of its industry uses are in the programming languages Erlang and Scala. Designed for concurrency, distribution, and scalability, actors are a part of Erlang itself. Scala (which stands for ``scalable language'') is a more recent programming language that runs on the Java platform and interoperates with all Java libraries, implementing actors using the multi-threaded capabilities of the Java Virtual Machine (JVM), while removing the complexity of thread communication for programmers \cite{scala} \cite{twitter-scal}. In fact, the popular social-networking site Twitter is using the Actor model for scalability by rewriting their backend infrastructure from Ruby on Rails to Scala \cite{twitter-scal}. Moreover, the Open Systems Laboratory, led by Gul Agha, aims to simplify the development of scalable parallel, distributed, and mobile computing systems, basing their research on the Actor model \cite{osl}.\newline
\indent Although interaction is not adequately encapsulated by the Turing machine and there are a multitude of conceptual frameworks (e.g. CSP, the $\pi$-calculus, Actor model, etc.) for dealing with interaction, there are still certain important aspects of computational practice that are not well expressed by these frameworks. The underlying conclusion is that these frameworks are intended to provide an \emph{abstract} framework for dealing with and reasoning about concurrency (and therefore abstract away from certain real world details). With the Actor model, for instance, timing details follow two general requirements. First, each actor always has the capability to process messages sent to it. Second, every message eventually arrives at its target, known as \emph{finite delay}\footnote{In fact, the type of computation that usually involves actors and distributed process models does not relate to single algorithmic computations, but can be viewed as functional compositions of such algorithms, where despite the latency of their result (namely, the second property of finite delay), these interactions are oracles, no different from regular function calls. Latency is a rather natural property of abstract descriptions of concurrent systems; finite but unbounded delay gives rise to unbounded nondeterminism. Although the Actor model exhibits unbounded nondeterminism (due to the underlying property of fairness), it must be made clear that this does \emph{not} mean that the Actor model can physically compute any functions outside the class of recursive functions.}. These requirements leave much ambiguity about the order of events in an actor's arrival ordering, and the resulting nondeterminism is called \emph{arrival nondeterminism}. The reason for the ambiguity in these requirements is simple, namely that: \newline
\begin{quote}
For a programming language semantics to specify completely the order in which events are to occur during multiprocessor execution of a concurrent program is generally impossible, since it would entail fixing myriad details such as the number and relative speeds of concurrent processors, the exact times of and delays occasioned by page faults and other interrupts, the timing of signals between processors and the manner in which they are arbitrated, and so on, down to the levels of time resolution at which quantum indeterminacy becomes important. Clinger \cite[pp. 56-7]{clinger} \newline
\end{quote}
However, in cyber-physical systems, there is a real need for theories that help in the design of these systems when not only interaction, but time plays a role; for example, at the device driver level, time is especially crucial. Time is not the only challenge that needs to be addressed. Networked Cyber-Physical Systems (NCPS) present many challenges not properly addressed by current distributed computing paradigms, such as limitations of computational, energy, and networking resources, as well as accounting for uncertainty, partial knowledge, and bad or stale information \cite{ncps}\cite{vision}. In relation to the usability constraint, NCPS are certainly ``usable'', where an emphasis is placed particularly on the sub-constraints of repeatability, physical constructibility, and reliability. Although work on cyber-physical systems deals with the engineering problem of integrating software with the real world, the problem lies deeper than that as ``it is becoming increasingly clear that more fundamental work is needed to address real-world challenges in a uniform and systematic way'' \cite{ncps}. The work of Stehr, et al. \cite{ncps}\cite{vision} aims to explore a new notion of software that behaves closer to a physical or biological system. This is motivated by the fact that in challenging environments, the physical world imposes severe limitations on how distributed algorithms can operate, and instead of viewing resource limitations as a bug, they view it as a feature -- and traditional models of distributed computing are too abstract to provide a proper foundation. The applications of NCPS are multifarious, providing complex, situation-aware services for traffic control, medical systems, automated laboratories, deep space exploration, emergency response, or social networking, to name a few. In fact, biological systems provide many examples of an effective NCPS, such as the human immune system, which as examples of desirable features has robustness, distributed knowledge, and generic and adaptive responses to events \cite[pg. 112]{vision}. There is much more that can be said about NCPS, so I refer the reader to Stehr, et al. \cite{vision} instead, where a long-term vision as well as research directions pertaining to NCPS are contained. The development and applications of NCPS offer an exciting venue for hypercomputation in terms of visiting the theoretical shortcomings of distributed computing paradigms and their adaptations to the physical world. While the abstract notions of relative computability have been around since the 1940s and were further greatly developed in the 1970s and 1980s, offering a mathematical framework for interactive and online computing that the \emph{a}-machine could not adequately provide, they are not sufficient to address the many challenges of real world computing, and this is not just an engineering problem. What is needed in distributed computing is a closer synergistic relationship between foundational concerns as well as a consideration of the real world details that are crucially relevant to computational practice.

\section{Conclusions}
\indent The hypercomputation movement, in its current form, is certainly not of interest to contemporary computer science practice. Since certain models of hypercomputation may hold considerable import to the theoretical aspects of computer science (as well as to debates on the foundations of physics, such as the PMH proposal), I am not saying that these models should be discarded altogether either. Rather, my interest stems from practical concerns regarding hypercomputability. Instead of asserting that non-computable functions, tasks, etc. are physically computable, a more reasonable extension to the classical paradigm should be focused on modeling computational processes that are not adequately encapsulated by the Turing machine (such as contemporary problems dealing with certain issues of interaction and tractability) -- the question of whether or not these processes allow one to compute the non-computable is not considered. \newline
\indent But why should my proposal be considered? Through a rigorous analysis of Turing computability and by adhering to rather general and reasonable epistemic contraints on physical computation, I demonstrated that Gandy's axiomatic principles for machine computation are sufficiently general given current physical theories. Considering that hypercomputation is focused on surpassing the limits on computation imposed by absolute computability (namely, computing more than the universal Turing machine), it is doubtful that even its goal would be of much interest to current computer science practice. Questions regarding issues of complexity as well as the notions of relative computability are examples of questions of prime interest in computational practice. In particular, my proposal contained issues for which Turing computability does not help, namely having to do with exploring alternative computational pradigms for effectively dealing with intractable problems in practice as well as building robust cyber-physical systems, and the scope of these issues are general enough for there to be ample opportunity for further research, yet mathematically precise enough for their importance to computational practice to be well-known and accepted. It is these aforementioned problems that a plausible extension to the classical paradigm should focus on in order to be practically viable.

\newpage 
\section*{Acknowledgements}
\indent I am grateful to Professor Solomon Feferman for his invaluable guidance throughout the process of this research. I also express my gratitude to Steven Ericsson-Zenith for discussions on interactive computing and the historical premises of Turing's work, and to Carolyn L. Talcott for discussions on the Actor model and cyber-physical systems. I further thank Martin Davis, Wilfried Sieg, and the two anonymous referees for their helpful comments on a draft of this article.
\newpage
\bibliographystyle{amsplain}

\end{document}